\documentclass[preprint,3p,times]{elsarticle}

\usepackage{amssymb}
\usepackage{amsmath}
\usepackage[all,cmtip]{xy}
\usepackage{amsthm}
\usepackage{color}
\usepackage{enumitem}
\usepackage{tikz}
\usepackage{appendix}
\usepackage[colorlinks=true]{hyperref}
\usepackage{hyperref}
 
\input diagxy
\theoremstyle{plain}
\newtheorem{thm}{Theorem}[section]
\newtheorem{lem}[thm]{Lemma}
\newtheorem{prop}[thm]{Proposition}
\newtheorem{cor}[thm]{Corollary}
\theoremstyle{definition}
\newtheorem{defn}[thm]{Definition}

\newtheorem{exmp}[thm]{Example}
\newtheorem{rem}[thm]{Remark}

\begin{document}

\newcommand{\ra}{\rightarrow}
\newcommand{\lra}{\longrightarrow}
\newcommand{\id}{\mathrm{id}}
\newcommand{\CL}{\mathcal{L}} 
\newcommand{\Ln}{\mathrm{Ln}}
\newcommand{\Exp}{\mathrm{Exp}}
\newcommand{\Max}{\mathrm{Max}}
 
\newcommand{\bv}{\bigvee}
\newcommand{\bw}{\bigwedge}
\newcommand{\Lra}{\Longrightarrow}
\newcommand{\sfM}{\mathsf{M}}
\newcommand{\Idm}{\mathbf{Idm}_{\otimes}}
\newcommand{\Lu}{\text{\L}}
\newcommand{\lam}{\lambda}
\newcommand{\de}{\delta}
\newcommand{\ga}{\gamma} 
\newcommand{\al}{\alpha}
\newcommand{\be}{\beta}
\newcommand{\Del}{\Delta}
\newcommand{\Delp}{\Delta^+}
\newcommand{\Nabp}{\nabla^+}
\newcommand{\CDp}{\mathcal{D}^+}
\newcommand{\vNabp}{\check{\nabla}^+}
\newcommand{\vep}{\varepsilon}
\newcommand{\Rop}{{[0,\infty]^{\rm op}}}
\newcommand{\vphi}{\varphi}
\newcommand{\SUP}{\mathbf{SUP}}
\newcommand{\INF}{\mathbf{INF}}

\begin{frontmatter}	
\title{Triangle functions generated by tensor products of quantales}

 \author[1]{Hongliang Lai}
\ead{hllai@scu.edu.cn}

\author[2]{Qingzhu Luo\corref{cor}}
\ead{qingzhuluo@cuit.edu.cn }

\cortext[cor]{Corresponding author.}
\address[1]{School of Mathematics, Sichuan University, Chengdu 610064, China}
\address[2]{College of Applied Mathematics, Chengdu University of Information Technology, Chengdu 610225, China}
\date{}
\begin{abstract}
This paper investigates triangle functions induced by tensor products of triangular norms and conorms.  For any left continuous t-norm $T$ on $[0,1]$ and any right continuous t-conorm $L$ on $[0,\infty]$, the tensor product $L\otimes T$ induces a triangle function on $\Delp$, giving rise to a partially ordered monoid structure on $(\Delta^+, L \otimes T)$.   The main results are as follows: (1) if $L$ is continuous, then $\tau_{T,L}$ is a triangle function on $\Delp$ if and only if $\tau_{T,L}=L\otimes T$, which in turn holds if and only if $L$ satisfies the property (LCS); (2) for $\CDp$, the set of all non-defective distance distribution functions, $(\CDp,L\otimes T)$ forms a submonoid of $(\Delp,L\otimes T)$ if and only if $L$ has no zero divisors; (3)for $\CDp_c$, the set of all continuous
distance distribution functions, if the t-norm $T$ is continuous, then $(\CDp_c,L\otimes T)$ is a subsemigroup of $(\Delp,L\otimes T)$ if and only if $L$ satisfies the property (LS). Furthermore, $(\CDp_c,L\otimes T)$ is an ideal of $(\CDp, L\otimes T)$ if and only if $L$ adheres
to the cancellation law.
\end{abstract}	

\end{frontmatter}

\section{Introduction}

In 1942, Menger introduced the probabilistic metric to quantify uncertainty in distances between elements of a space \cite{Menger1942}. Rather than representing distance as a deterministic real value, he proposed characterizing it through a distance distribution function (d.d.f.) $F_{pq}$, which treats the distance between points $p$ and $q$ as a random variable with an associated probability distribution. This framework was later advanced by \v{S}erstnev in 1962, who formalized the notion of a triangle function defined over the set of all distance distribution functions. By introducing a probabilistic analogue of the triangle inequality, \v{S}erstnev extended classical metric axioms to probabilistic settings, thereby contributing to the rigorous development of probabilistic metric spaces \cite{Serstnev1962}.


Let $\Delp$ denote the set of all distance distribution functions. Under the pointwise order, $\Delta^+$ forms a complete lattice. Informally, a triangle function is a binary operation on $\Delp$ that plays a role analogous to the addition of real numbers, making it a fundamental concept in the theory of probabilistic metric spaces. In \cite{Schweizer1983}, Schweizer and Sklar compiled various methods for constructing triangle functions and presented several distinct classes of triangle functions. Among these, an important family of triangle functions, denoted $\tau_{T,L}$, is generated by a left continuous triangular norm $T$ on $[0,1]$ and a continuous and commutative monoid operation $L$ on $[0,\infty]$ with unit element $0$ and satisfying the property
\[x<y, x'<y'\implies L(x,y)<L(x',y'). \leqno{(LS)}\]
Specifically, for any distance distribution functions $f$ and $g$, and for any $x\in [0,\infty]$, we define
$$\tau_{T,L}(f, g)(x) = \sup \left\{ T(f(u), g(v)) \mid L(u, v) = x \right\},$$
then $\tau_{T,L}$ is a triangle function on $\Delta^+$ (see Theorem 7.2.4 in \cite{Schweizer1983}).

Within the framework of quantale theory, a novel method for constructing triangle functions on $\Delp$ is established. It is observed that in the category $\mathbf{SUP}$ of complete lattices and supremum-preserving maps, $\Delp$ can be identified with the tensor product of the complete lattices $[0,1]$ and $[0,\infty]^{\mathrm{op}}$. Moreover, let $T$ be a left continuous triangular norm on $[0,1]$, and let $L$ be a right continuous commutative monoid operation on $[0,\infty]$ with identity element $0$. Then $([0,1],T,1)$ and $([0,\infty]^{\mathrm{op}}, L, 0)$ form quantales. The tensor product of these quantales yields a complete lattice $\Delp$, equipped with a monoid operation $L \otimes T$, which is a triangle function on $\Delp$.

The two methods for constructing triangle functions described above appear distinctly different, yet they can yield the same triangle function under certain conditions. For instance, when $L$ coincides with the addition operation on the real numbers, Gutiérrez García et al. observed that $\tau_{T,L}=L\otimes T$ \cite{Eklund2018, Gutierrez2017, Gutierrez2019}. This observation naturally leads to the following question: under what general conditions does $\tau_{T,L} = L \otimes T$ hold? We prove that for any left continuous t-norm $T$ on $[0,1]$ and any continuous t-conorm $L$ on $[0,\infty]$, the equality $\tau_{T,L} = L \otimes T$ holds if and only if  $L$ satisfies the  property (LCS).  This result not only improves upon Theorem 7.2.4 in \cite{Schweizer1983}, but also partially answers the first question posed in Problem 7.9.6 of \cite{Schweizer1983}.

For a triangle function $\tau$ on $\Delta^+$, the study of subsets of $\Delta^+$ that are closed under the operation $\tau$ presents an interesting problem.  In Problem 7.9.3 (i) of \cite{Schweizer1983}, the following open question is included:
``Given a (continuous) triangle function $\tau$, determine the subsets of $\Delta^+$ closed under $\tau$. In particular, is $\mathcal{D}^+$ closed under $\tau$?'' There are few sufficient conditions for $\CDp$ being closed  listed in \cite[Theorem 12.1]{Saminger2008}. As for the subset $\CDp_c$ of continuous distance distribution functions, it is shown in \cite{Moynihan1975} that $\tau_{T,L}(f,g)$ is continuous whenever $f,g$ are in $\Delp$ and either one is continuous, where $T$ is a continuous t-norm and $L$ is the addition operation.

In this paper, we focus on triangle functions of the form $L \otimes T$ and examine whether the following subsets of $\Delta^+$ are closed under this class of triangle functions:

\begin{itemize}
\item $\mathcal{D}^+$, the set of all non-defective distance distribution functions;
\item $\mathcal{D}_0^+$, the set of all distance distribution functions that are continuous on $]0, \infty]$;
\item $\mathcal{D}_c^+$, the set of all distance distribution functions that are continuous on $[0, \infty]$;
\end{itemize}

We demonstrate that the closedness of the aforementioned subsets  under the operation $L\otimes T$ generally requires the t-norm $T$ and the t-conorm $L$ to satisfy certain additional properties. The main results are as follows:

\begin{enumerate}
\item[(i)] $\CDp$ is closed under $L\otimes T$ if and only if $L$ has no zero divisors.
\item[(ii)] $\CDp_0$ is closed under $L\otimes T$ if and only if $T$ is continuous and $L$ satisfies the property  $(LS)$.
\item[(iii)] If $T$ is continuous, then $\CDp_c$ is closed under $L\otimes T$ if and only if $L$ satisfies the property $(LS)$.
\item[(iv)] $L\otimes T(f,g)\in\CDp_c$
  for all $f\in\CDp$ and $g\in\CDp_c$
if and only if $T$ is continuous and $L$ satisfies the cancellation law.
\end{enumerate}

 The paper is structured as follows: In Section 2, we review fundamental concepts related to distance distribution functions and tensor products of complete lattices. In Section 3, we analyze the relationship between the operation $\tau_{T,L}$ and the tensor product operation $L\otimes T$. In Section 4, we investigate the closedness of the subsets $\CDp$, $\CDp_0$ and $\CDp_{c}$ under the operation $L\otimes T$.

\section{Distance distribution functions}

Let $[0,\infty]$ be the extended non-negative reals, i.e., the set of non-negative real numbers with the symbol $\infty$ added, and ordered by $0\leq x <\infty$ for all  $x\in[0,\infty[$. 

Let $[a,b]$ be a closed interval contained in $[0,\infty]$. The supremum of a set $\{x_i:i\in I\}\subseteq[a,b]$ is denoted by $\sup\limits_{i\in I}x_i$ and the infimum is denoted by $\inf\limits_{i\in I}x_i$. Particularly, if the index set $I$ is empty, $\sup_{i\in I}x_i=a$ and $\inf\limits_{i\in I}x_i=b$. 

Let $[a,b]$ and $[c,d]$ be closed intervals contained in $[0,\infty]$. An increasing function 
$f:[a,b]\lra [c,d]$ is  supremum-preserving   if it preserves arbitrary suprema, that is, 
 $$f(\sup_{i\in I} x_i)=\sup_{i\in I}f(x_i)$$ 
 for all family $\{x_i: i\in I\}$ in $[a,b]$. Particularly, if the index set $I$ is the empty set, one has that $f(a)=c$. Dually,  $f$ is infimum-preserving if it preserves all infima, that is, $$f(\inf_{i\in I} x_i)=\inf_{i\in I}f(x_i)$$ for all family $\{x_i: i\in I\}$ in $[a,b]$. If the index set $I$ is the empty set, one has that $f(b)=d$.

 Given an increasing function $f:[a,b]\lra[c,d]$, define
 \[f^{-}(x)=\sup_{y<x} f(y), \quad f^{+}(x)=\inf_{y>x}f(y),\]
 then $f^-$ is a supremum-preserving map  and $f^+$ is an infimum-preserving map. In fact, $f^-$ is the largest supremum-preserving map  below $f$ and   $f^+$ is the smallest infimum-preserving map    above $f$. Hence,  $f=f^-$ if and only if $f$ is a supremum-preserving map , and $f=f^+$ if and only if  $f$ is an infimum-preserving map. Furthermore, one can easily see that, $f^{-}\leq g\iff f\leq g^{+}$ for all increasing functions $f,g:[a,b]\lra [c,d]$,  $(f^+)^-=f$  for all supremum-preserving function $f:[a,b]\lra[c,d]$, and 
 $(f^-)^+=f$ for all infimum-preserving function $f:[a,b]\lra[c,d]$.

\begin{defn}\cite{Schweizer1983}
A \emph{distance distribution function} (briefly, a \emph{d.d.f.}) is an increasing function $f:[0,\infty]\lra [0,1]$ 
satisfies $f(0)=0$ and $f(\infty)=1$, and is left continuous on the open interval $]0,\infty[$. The set of all d.d.f.'s is denoted by $\Delp$.
\end{defn}

\begin{rem}\label{Delp}
If $f:[0,\infty]\lra[0,1]$ is a d.d.f., then 
$$f^-(x)=\begin{cases} f(x) & x<\infty,\\
                       \sup\limits_{t<\infty}f(t)    & x=\infty, 
        \end{cases}
$$ and $f^-$ is a supremum-preserving.
Conversely, given a supremum-preserving map  $g:[0,\infty]\lra[0,1]$,    define a function 
$$\widetilde{g}(x)=\begin{cases}g(x)& x<\infty,\\
                            1 &   x=\infty.
                \end{cases}$$

This establishes a bijective correspondence:
\begin{itemize}
\item $\widetilde{f^-}=f$ for all $f\in\Delp$,
\item $(\widetilde{g})^-={g}$  for all supremum-preserving map  $g:[0,\infty]\lra[0,1]$. 
\end{itemize}
Thus, $\Delp$ is isomorphic to the set of all supremum-preserving maps from $[0,\infty]$ to $[0,1]$. It is worth noting that some literature uses $\Delp$ exclusively for this class of maps \cite{Eklund2018, Gutierrez2017, Gutierrez2019}. However, in this paper, we adhere to the notation in \cite{Saminger2008, Saminger2010,Schweizer1983} to facilitate direct comparison between the tensor product $L\otimes T$ and the triangle function
$\tau_{T}^{L}$ established in existing works.
\end{rem}

\begin{defn}\cite{Schweizer1983}  Given a d.d.f. $f:[0,\infty]\lra[0,1]$, the largest quasi-inverse of $f$ is a function 
$$f^\vee:[0,1]\lra[0,\infty],\quad f^\vee(y)=\inf\{x\in[0,\infty]\mid f(x)> y\}.$$
All largest quasi-inverses of d.d.f.'s  form a set, denoted by 
$$\Nabp=\{f^\vee\mid f\in\Delp\}.$$
\end{defn}


\begin{prop} {\rm(\cite{Schweizer1983})}\label{pro adj}
Let $f:[0,\infty]\lra[0,1]$ be a d.d.f..    
 Then $f^\vee(y)=\sup\{x\in[0,\infty]\mid f(x)\leq y\}$ and $f^-(x)\leq y\iff x\leq f^\vee(y)$ 
for all $x\in[0,\infty]$ and all $y\in[0,1]$. 
\end{prop}


\begin{rem} For d.d.f.'s $f$ and $g$, one has that 
$$f\leq g\iff f^-\leq g^-\iff f^\vee\geq g^\vee.$$
Therefore, $\Nabp$ is dually isomorphic to $\Delp$. 
\end{rem}

\begin{prop}{\rm(\cite{Schweizer1983})}\label{cont vs stri incr} Let $f$ be a d.d.f.. Then 
\begin{enumerate}
    \item[(i)] $f$ is continuous on a closed interval $[a,b]$ if and only if $f^\vee$ is increasing strictly on the interval $[f(a),f(b)]$. Consequently, $f$ is continuous on an open interval $]a,b[$ if and only if  $f^\vee$ is increasing strictly on the interval $]f(a^+),f(b^-)[$. 
    \item[(ii)] $f$ is increasing strictly on $[a,b]$ if and only if $f^\vee$ is continuous on the interval $[f(a),f(b)]$.
\end{enumerate}
\end{prop}






\begin{exmp} For all $r\in[0,\infty]$ and all $p\in[0,1]$, define a step function 
$$\vep(r,p)(x)=\begin{cases} 1  & x=\infty,\\
                             p  & r<x<\infty,\\
                             0  & 0\leq x\leq r.
                \end{cases}$$
It is easy to see that if $r=\infty$ or $p=0$, then
$$\vep(r,p)(x)=\begin{cases}1 &  x=\infty,\\
                            0 &  0\leq x<\infty.
            \end{cases}
$$                        

For each step function $\vep(r,p)$, we calculate the following functions explicitly:

$$\vep^-(r,p)(x)=\begin{cases} p& r<x\leq\infty,\\
                              0& 0\leq x\leq r;
                \end{cases}
\quad 
\vep^\vee(r,p)(y)=\begin{cases}\infty& p\leq y\leq 1,\\
                                      r& 0\leq y<p.
                \end{cases}$$
  

\end{exmp}

Equipped with the pointwise order, $\Delp$ is a complete lattice, in which   $\vep(0,1)$ is the top element and $\vep(\infty,0)$ is the bottom element. The supremum of a non-empty family $\{f_i:i\in I\}$ is calculated pointwisely, that is, 
$$(\sup_{i\in I}f_i)(x)=\sup_{i\in I}f_i(x)$$ for a family of d.d.f.'s $\{f_i, i\in I\}$ in $\Delp$.  

\begin{prop}\label{equ of ddf}
Let $f$ be a d.d.f.. Then it holds that 
\begin{enumerate}
\item[{\rm(1)}] $f=\sup\limits_{r\in[0,\infty]}\vep(r, f(r))=\sup\limits_{p\in[0,1]}\vep(f^\vee(p),p)$;
\item[{\rm(2)}] $f^\vee=\inf\limits_{r\in[0,\infty]}\vep^\vee(r,f(r))
=\inf\limits_{p\in[0,1]}\vep^\vee(f^\vee(p),p)$.
\end{enumerate}
\end{prop}
\begin{proof}
$(1)$ The equality holds trivially at the point $0$ and $\infty$. For all $x\in ]0,\infty[$, $$\sup\limits_{r\in[0,\infty]}\vep(r, f(r))(x)=\sup\limits_{r<x}f(r)=f(x)$$ since $f$ is left continuous on $]0,\infty[$.
By Proposition \ref{pro adj} $(1)$, for all $x\in]0,\infty[, y\in[0,1]$, \[f(x)\leq y\iff x\leq f^\vee(y),\]
which implies $f(x)=\inf\{y\mid f^\vee(y)\geq x\}$. Hence $$\sup\limits_{p\in[0,1]}\vep(f^\vee(p),p)(x)=\sup\{p\mid f^\vee(p)<x\}=\inf\{p\mid f^\vee(p)\geq x\}=f(x)$$  for all $x\in]0,\infty[$.

    $(2)$ Since $\Delp$ and $\Nabp$ are dually isomorphic, it follows that 
    $$(\sup\limits_{i\in I}f_i)^\vee=\inf\limits_{i\in I}{f_i}^\vee$$ for any family of d.d.f.'s $\{f_i,i\in I\}$ in $\Delp$, and then the equality holds dually. 
\end{proof}

A complete lattice $M$ is a partially ordered set which admits suprema of arbitrary subsets of $M$. A morphism $f:M\lra N$ between complete lattices is a supremum-preserving map. The category of all complete lattices and supremum-preserving maps is denoted by $\mathbf{SUP}$.  

Let $M_1$, $M_2$ and $M_3$ be complete lattices.   A map $\vphi:M_1\times M_2\lra M_3$ is a bi-morphism  if it preserves suprema in each place separately, that is,
 $$\vphi(\sup_{i\in I}x_i,y)=\sup_{i\in I}\vphi(x_i,y),\quad 
 \vphi(x,\sup_{i\in I}y_i)=\sup_{i\in I}\vphi(x,y_i).$$

\begin{defn}\cite{Banaschewski1976,Joyal1984}  The tensor product $M_1\otimes M_2$ of complete lattices $M_1,M_2$ in the category $\SUP$ is the codomain of the universal bi-morphism  $$\eta: M_1\times M_2\lra M_1\otimes M_2.$$ That is, for each bi-morphism  $\vphi:M_1\times M_2\lra M_3$, there is a unique supremum-preserving map  $\widehat{\vphi}:M_1\otimes M_2\lra M_3$ such that $\vphi=\widehat{\vphi}\circ\eta$. 
\begin{center}
\begin{tikzpicture}
\draw[->] (0.8,2)--(2,2) node[right] {$M_3$};
\draw[->] (0,1.6)--(0,0) node[below] {$M_1\otimes M_2$};
\draw[->] (0.8,0)--(2,1.6);
\node at(0,2){$M_1\times M_2$};
\node at(0,.8)[left]{$\eta$};
\node at(1.4,2)[above]{$\varphi$};
\node at(1.4,.8)[right]{$\widehat{\varphi}$};
\end{tikzpicture}
\end{center}
\end{defn}

\begin{rem}
The tensor product of complete lattices $M_1$ and $M_2$
is unique up to isomorphism, implying multiple different representational forms for $M_1\otimes M_2$. Let $M_1=\Rop$ (the interval $[0,\infty]$ with the reversed order, hence, $\Rop$ has $0$ as its top element and $\infty$ as its bottom element.) and $M_2=[0,1]$.  Two distinct representations of $\Rop\otimes[0,1]$ are: The complete lattice of all infimum-preserving maps $f:[0,\infty]\lra[0,1]$ and 
the complete lattice of all supremum-preserving maps $g:[0,\infty]\lra[0,1]$. See Example 2.1.10 in \cite{Eklund2018} and Subsection 3.3 in \cite{Gutierrez2017} for details.

By adjusting the value $f(\infty)$, a supremum-preserving map $f:[0,\infty]\lra[0,1]$ can be bijectively converted into a d.d.f. and vice versa (see Remark \ref{Delp}). This establishes that $\Delp$ is another canonical representation of $\Rop\otimes[0,1]$.
Clearly, the  universal  bi-morphism   is given by 
$$\vep: \Rop\times[0,1]\lra\Delp,\quad (r,p)\mapsto \vep(r,p).$$    
Given any bi-morphism  $\vphi:\Rop\times[0,1]\lra M$   with the codomain $M$, the map 
$$\widehat{\vphi}: \Delp\lra M, \quad
  \widehat{\vphi}(f)=\bigvee_{r\in[0,\infty]}\vphi(r,f(r))$$
is the unique supremum-preserving map  such that $\vphi=\widehat{\vphi}\circ\vep$.
\end{rem}
  

\section{Triangle functions generated by t-norms on $[0,1]$ and t-conorms on $[0,\infty]$}

 \begin{defn} \cite{Schweizer1983}  A \emph{triangle function} is an associative  binary operation $\tau$ on $\Delp$ that is commutative, increasing in each place, and has $\vep(0,1)$ as identity.\end{defn}

 \begin{defn}\cite{Durante2005,Durante2006} 
 A \emph{semicopula} is a function $S:[0,1]\times[0,1]\lra [0, 1]$  that satisfies the
following two conditions:
\begin{itemize}
\item[(i)] $S$ is increasing in each place, 
\item[(ii)] for all $x\in[0,1]$, $ S(x,1)=S(1,x)=x$. 
\end{itemize}
The set of all semicopulas is denoted by $\mathcal{S}$. 
 \end{defn}

\begin{defn} \cite{Frank1979,Saminger2008,Saminger2010,Schweizer1983} The class $\CL$ is the set of all binary operations $L$ on $[0,\infty]$ such that:
\begin{enumerate}
    \item [(\rm L1)] $L:[0,\infty]\times[0,\infty]\lra[0,\infty]$ is onto;
    \item [(\rm L2)] $L$ is increasing in each place;
    \item [(\rm L3)] $L$ is continuous on $[0,\infty]\times{[0,\infty]}$, except possibly at the points $(0,\infty)$ and $(\infty,0)$.
\end{enumerate}
\end{defn}

On occasion, we impose additional properties on a binary operation $L$ on $[0,\infty]$. These include the following:
\begin{itemize}
\item $L$ has the identity $0$, if  
\[\forall x\in[0,\infty], L(0,x)=L(x,0)=x. \leqno{(L0)}\]
\item  $L$ satisfies the strictly increasing property,   if 
$$\forall x,x',y,y'\in[0,\infty], x<x', y<y'\implies L(x,y)<L(x',y').\leqno(LS)$$
\item $L$ satisfies the conditionally strictly increasing property, if 
$$\forall x,x',y,y'\in[0,\infty], x<x', y<y', L(x',y')<\infty\implies L(x,y)<L(x',y').\leqno(LCS)$$
\end{itemize}

\begin{defn}{\cite{Schweizer1983}} Let $T$ be a semicopula and $L$ an element in the class $\CL$. For any $f,g$ in $\Delp$, define a function $\tau_{T,L}(f,g):[0,\infty]\lra[0,1]$ by 
$$\tau_{T,L}(f,g)(x)=\sup_{L(s,t)=x}T(f(s),g(t)).$$
\end{defn}

Generally speaking, $\tau_{T,L}$
  is not necessarily a binary operation on $\Delp$, since there exist distance distribution functions $f$ and $g$ for which $\tau_{T,L}(f,g)$ fails to be a distance distribution function.
 In order for $\tau_{T,L}$ to form a triangle function on $\Delp$, additional properties must be imposed on $T$ and $L$.

\begin{defn} \cite{Alsina2006,Klement2000,Schweizer1958,Schweizer1960,Schweizer1961,Schweizer1983} A \emph{triangular norm} (shortly, a t-norm) is an associative binary operation $T$ on $[0,1]$ that  is   commutative, increasing in each place, and has $1$ as identity. A t-norm $T$ is \emph{left continuous} (or  \emph{continuous}) if $T: [0,1]^2\lra [0,1]$ is a left continuous (or   continuous) function with respect to the usual topology on $[0,1]$.
\end{defn}

The definition of a triangular conorm on the interval $[0,1]$ can be adapted to a bounded partially ordered set by making slight modifications \cite{DeCooman1994,Drossos1999}. In this paper, we defined the triangular conorm on the interval $[0,\infty]$ as follows:

\begin{defn} A \emph{triangular conorm} (shortly, a t-conorm) is an associative binary operation $L$ on $[0,\infty]$ that  is   commutative, increasing in each place, and has $0$ as identity. A t-conorm $L$ is \emph{right continuous} (or  \emph{continuous}) if $L: [0,\infty]\times[0,\infty]\lra [0,\infty]$ is a right continuous (or   continuous) function with respect to the usual topology on $[0,\infty]$.
\end{defn}

\begin{prop}{\rm(\cite{Schweizer1983}, Theorem 7.2.4)}\label{tauTL is tria func}
    Let $T$ be a left continuous t-norm. Let $L$ be a commutative and associative operation in $\CL$ which additionally satisfies properties (LS) and (L0). Then the operation $\tau_{T,L}$ is a triangle function on $\Delp$.
\end{prop}

\begin{rem}
A commutative and associative operation $L$ in $\CL$ satisfies the property (L0) if and only if it is a continuous t-conorm on $[0,\infty]$. In fact, the ``if'' part is obvious. To check the ``only if'' part, firstly, notice that the property (L0) yields that $L$ has the identity element $0$.  Secondly, we claim that $L$ is continuous at the points $(0,\infty)$ and $(\infty,0)$, hence, it is a continuous binary operation on $[0,\infty]$. In fact, it holds that $$\lim_{(x,y)\ra(\infty,0)}L(x,y)=L(\infty,0)=\infty=L(0,\infty)=\lim_{(x,y)\ra(0,\infty)}L(x,y),$$ since $L(x,y)\geq L(0,y)=y$ and $L(x,y)\geq L(x,0)=x$. 
\end{rem}

Therefore, the aforementioned proposition can be reformulated as follows:
\begin{prop}\label{tauTL is tria func refo}
Let $T$ be a left continuous t-norm on $[0,1]$ and $L$ be a continuous t-conorm on $[0,\infty]$ satisfying the property (LS). Then $\tau_{T,L}$ is a triangle function on $\Delp$.
\end{prop}



There are three basic t-norms on the unit interval $[0,1]$ in the following:
\begin{itemize}
\item The G\"{o}del t-norm: $M(x,y)=\min\{x,y\}$;
\item The product t-norm: $\Pi(x,y)=x\cdot y$; 
\item The {\L}ukasiewicz t-norm: $W(x,y) =\max\{x+y-1, 0\}$.
\end{itemize}

Supplement  $-\ln(0)=\infty$ and $\mathrm{e}^{-\infty}=0$,
then the functions $x\mapsto \mathrm{e}^{-x}$ is the inverse of   $y\mapsto -\ln y$ between the intervals $[0,\infty]$ and $[0,1]$. Through this pair of dually isomorphisms, each t-norm is transposed to a t-conorm
$$T^*(x,y)=-\ln T(\mathrm{e}^{-x},\mathrm{e}^{-y}).$$
For example, 
\begin{itemize}
\item $M^*(x,y)=-\ln M(\mathrm{e}^{-x},\mathrm{e}^{-y})=\max\{x,y\}$,
\item $\Pi^*(x,y)=-\ln \Pi(\mathrm{e}^{-x},\mathrm{e}^{-y})=x+y$.
\end{itemize}

Let $L$ be a t-conorm on $[0,\infty]$. It satisfies the cancellation law if 
$$  \forall x,y\in[0,\infty], L(x,y)=L(x,z)\Lra x=\infty\text{ or } y=z. \leqno(CL)$$
It has no zero divisor if 
$$L(x,y)=\infty\implies x=\infty \text{ or } y=\infty.\leqno(NZD)$$
It is easily seen that, for each t-conorm $L$ on $[0,\infty]$, it holds that  
$$(CL)\implies(LS)\iff (LCS)\text{ and }(NZD).$$

\begin{exmp} 
\begin{enumerate}
\item The t-conorm  $W^*$  satisfies the condition $(LCS)$ but neither the condition $(LS)$ nor the condition $(NZD)$.  The t-conorm $M^*$ satisfies the condition $(LS)$ but not the cancellation law $(CL)$. The addition operation $+$ is a t-conorm satisfying the cancellation law. 

\item Define a continuous t-norm $T$ on $[0,1]$ by
$$T(x,y)=\begin{cases} 2x\cdot y & x,y\in[0,0.5],\\
                             \max\{x+y-1, 0.5\} & x,y\in[0.5,1],\\
        \min\{x,y\} & \text{otherwise}.
                \end{cases}$$
Then the continuous t-conorm $T^*$ has no zero divisor but does not satisfy the condition $(LCS)$.

\item There is a right continuous but not continuous t-conorm satisfying the cancellation law, which is isomorphic to the t-norm given in Example 2.5.2 in  \cite{Alsina2006}. Through the bijection $$\varphi:[0.5,1]\lra{[0,1]}, x\mapsto 2x-1,$$ one can define a binary operation $F$ on $[0.5,1]$ by 
$$F(p,q)=2(p\cdot q)-(p+q)+1,$$
which is isomorphic to the product t-norm $\Pi$ on $[0,1]$. Define a t-norm $T$ on $[0,1]$ by
$$ T(x,y)=\frac{F(2^m\cdot x,2^n\cdot y)}{2^{m+n}}, 
\quad x\in\bigg]\frac{1}{2^{m+1}},
\frac{1}{2^m}\bigg], y\in\bigg]\frac{1}{2^{n+1}},\frac{1}{2^n}\bigg]$$
for all natural numbers $m,n=0,1,2,\cdots$. One can see that $T$ is not continuous but left continuous,  and it satisfies the cancellation law in the sense of that
$$T(x,y)=T(x,z)\implies x=0\text{ or } y=z.$$

  Clearly, the t-conorm $T^*$ is right continuous but not continuous, and  satisfies the 
cancellation law $(CL)$ since $T$ does. 
\end{enumerate}
\end{exmp}

An \emph{integral and commutative quantale} \cite{Rosenthal1990} is a commutative monoid $(Q,\&)$ such that $Q$ is a complete lattice, and $\&$ has the top element of $Q$ as identity and preserving all suprema in each place. If $T$ is a t-norm, then $([0,1],T)$ is an integral and commutative quantale if and only if $T$ is left continuous. Similarly, if $L$ is a t-conorm on $[0,\infty]$, then $([0,\infty]^{\rm op},L)$ is an integral and commutative quantale if and only if $L$ is right continuous. 

 Recall from \cite{Eklund2018,Rosenthal1996} that we can construct their tensor product $(\Delp,L\otimes T)$, which is also an integral and commutative quantale, hence   the binary operation   $L\otimes T$ is a triangle function on $\Delp$. Precisely, $L\otimes T$ is a bi-morphism satisfying that 
 $$L\otimes T (\vep(r,p),\vep(s,q))=\vep(L(r,s),T(p,q))$$ 
 for all   step functions $\vep(r,p)$ and $\vep(s,q)$ in $\Delp$.
Thus, the value $L\otimes T(f,g)$ for all d.d.f.'s $f$ and $g$ is determined by 
$$L\otimes T(f,g)=\sup_{r,s\in[0,\infty]} L\otimes T(\vep(r,f(r)), \vep(s,g(s)))=\sup_{r,s\in[0,\infty]}\vep(L(r,s),T(f(r),g(s))).$$  

\begin{prop}\label{L tens T}
Let $T$ be a left continuous t-norm on $[0,1]$ and $L$ be a right continuous t-conorm on $[0,\infty]$. Then 
$$L\otimes T(f,g)(x)=\begin{cases}
 1 & x=\infty,\\
 \sup\limits_{L(r,s)<x}T(f(r),g(s)) & 0\leq x<\infty.
 \end{cases} $$ for all $f,g\in\Delp$.
 \end{prop}
\begin{proof} 
For all $x\in[0,\infty]$, we calculate explicitly that 
 \begin{align*}
 L\otimes T(f,g)(x)
 &=\sup_{r,s\in[0,\infty]}\vep(L(r,s),T(f(r),g(s)))(x)\\
 &=\sup_{L(r,s)<x}\vep(L(r,s),T(f(r),g(s)))(x)\\
 &=\begin{cases}
 1 & x=\infty,\\
 \sup\limits_{L(r,s)<x}T(f(r),g(s)) & 0\leq x<\infty.
 \end{cases} 
 \end{align*}
 \end{proof}

 \begin{rem} 
 Notice that for  d.d.f.'s $f$ and $g$ which are not left continuous at $\infty$,  the value $L\otimes T(f,g)(\infty)$ differs from that in \cite{Eklund2018, Gutierrez2017, Gutierrez2019}. 
 \end{rem}

\begin{prop}\label{Tens in Nabp}
Let $T$ be a left continuous t-norm on $[0,1]$ and $L$ be a right continuous t-conorm on $[0,\infty]$. Then it holds that
$$(L\otimes T(f,g))^\vee(y)=\inf_{T(p,q)>y}L(f^\vee(p),g^\vee(q))$$
for all $f,g\in\Delp$ and all $y\in[0,1]$. 
\end{prop}
\begin{proof}
For all $y\in [0,1]$, by Proposition \ref{equ of ddf}, one can calculate as:
\begin{align*}
 (L\otimes T(f,g))^\vee(y)&=[L\otimes T(\sup_{p\in[0,1]}\vep(f^\vee(p),p),\sup_{q\in[0,1]}\vep(g^\vee(q),q))]^\vee(y)\\
 &=[\sup_{p,q\in[0,1]}L\otimes T(\vep(f^\vee(p),p),\vep(g^\vee(q),q))]^\vee(y)\\
 &=[\sup_{p,q\in[0,1]}\vep(L(f^\vee(p),g^\vee(q)),T(p,q))]^\vee(y)\\
 &=\inf_{p,q\in[0,1]}\vep^\vee(L(f^\vee(p),g^\vee(q)),T(p,q))(y)\\
 &=\inf_{T(p,q)>y}\vep^\vee(L(f^\vee(p),g^\vee(q)),T(p,q))(y)\\
 &=\inf\limits_{T(p,q)>y}L(f^\vee(p),g^\vee(q)).
 \end{align*}
\end{proof}

\begin{exmp} Let $T$ be a left continuous t-norm on $[0,1]$ and $L$ be a right continuous t-conorm on $[0,\infty]$.      Then, for all $x\in[0,\infty[$, 
\begin{itemize}
\item[(i)] $M^*\otimes T(f,g)(x)=\sup\limits_{\max\{r,s\}<x}T(f(r),g(s))=\sup\limits_{r,s<x}T(f(r),g(s))=T(f(x),g(x))$, 
\item[(ii)] $(L\otimes M(f,g))^\vee(y)=\inf\limits_{\min\{p,q\}>y}L(f^\vee(p),g^\vee(q))=\inf\limits_{p,q>y}L(f^\vee(p),g^\vee(q))=L(f^\vee(y),g^\vee(y))$.
\end{itemize}
\end{exmp}


By making a slight modification to the definition of $\tau_{T,L}$, we define a new function, denoted $\tau_T^L$, as follows:
\begin{defn} Let $T$ be a left continuous t-norm on $[0,1]$ and $L$ be a right continuous t-conorm on $[0,\infty]$. For all $f,g\in\Delp$, define a function $\tau_T^L(f,g):[0,\infty]\lra[0,\infty]$ by 
$$\tau_{T}^{L}(f,g)(x)=\sup\limits_{L(r,s)\leq x}T(f(r),g(s)). $$
 \end{defn}

\begin{prop}\label{tau_T^L vs tau_T,L}
Let $T$ be a left continuous t-norm on $[0,1]$ and $L$ be a continuous t-conorm on $[0,\infty]$.  Then for all $f,g\in\Delp$, $\tau_{T}^{L}(f,g)=\tau_{T,L}(f,g)$. 
\end{prop}
  \begin{proof} 
 It suffices to show that  $$\sup\limits_{L(r,s)\leq x}T(f(r),g(s))\leq\sup\limits_{L(r,s)=x}T(f(r),g(s))$$ 
 for all $x\in[0,\infty[$. Let $x_0\in [0,\infty[$ and $r_0,s_0$ with $L(r_0,s_0)\leq x_0$. Since $L$ is continuous, there exists a point $t_0\in [s_0,\infty]$ such that $x_0=L(r_0,t_0)$, then $$T(f(r_0),g(s_0))\leq T(f(r_0),g(t_0))\leq\sup\limits_{L(r,s)=x_0}T(f(r),g(s)).$$ 
 Hence $$\sup\limits_{L(r,s)\leq x}T(f(r),g(s))\leq\sup\limits_{L(r,s)=x}T(f(r),g(s))$$
 as desired. 
 \end{proof}
 
\begin{thm}\label{tens vs tau} Let $T$ be a left continuous t-norm on $[0,1]$ and $L$ be a right continuous t-conorm on $[0,\infty]$. The following statements are equivalent:
 \begin{enumerate}
\item[{\rm (1)}]  for all d.d.f.'s $f$ and $g$, $\tau_{T}^{L}(f,g)=L\otimes T(f,g)$;

\item[{\rm (2)}]  for all d.d.f.'s $f$ and $g$, $\tau_{T}^{L}(f,g)$ is also a d.d.f.;

\item[{\rm (3)}] the t-conorm $L$ on $[0,\infty]$ fulfills the property $(LCS)$. 
\end{enumerate}
\end{thm}
\begin{proof}
(1)$\iff$ (2): On one hand, if $\tau_{T}^{L}(f,g)=L\otimes T(f,g)$, then $\tau_{T}^{L}(f,g)$ is clearly a d.d.f. since $L\otimes T$ is a binary operation on $\Delp$. On the other hand, if $\tau_{T}^{L}(f,g)$ is a d.d.f., then it is left continuous on the open interval $]0,\infty[$. Thus, for all $x\in]0,\infty[$, 
\begin{align*}
\tau_{T}^{L}(f,g)(x)&=\sup\limits_{z<x}\tau_{T}^{L}(f,g)(z)\\
&=\sup\limits_{z<x}\sup\limits_{L(r,s)\leq z}T(f(r),g(s))\\
&=\sup\limits_{L(r,s)<x}T(f(r),g(s))\\
&=L\otimes T(f,g)(x).
\end{align*}
Additionally, it always holds that $$\tau_{T}^{L}(f,g)(0)=L\otimes T(f,g)(0)=0.$$ 
and 
$$\tau_{T}^{L}(f,g)(\infty)=L\otimes T(f,g)(\infty)=1.$$ Therefore, $\tau_{T}^{L}(f,g)=L\otimes T(f,g)$ as desired. 

(2)$\implies$(3): If $L$ does not satisfy the condition $(LCS)$, then there are $r,r',s,s'$ in $[0,\infty]$ such that $r<r'$,  $s<s'$ and $L(r,s)<\infty$, but $L(r,s)=L(r',s')$. We consider the function $\tau_{T}^{L}(\vep(r,1),\vep(s,1))$ in the below.

Let $x_0=L(r,s)=L(r',s')$, then $x_0\geq\max\{r',s'\}> 0$.   We obtain that 
\begin{align*}
\tau_{T}^{L}(\vep(r,1),\vep(s,1))(x_0)&=\sup_{L(y,z)\leq x_0}T(\vep(r,1)(y),\vep(s,1)(z))\\
&\geq T(\vep(r,1)(r'),\vep(s,1)(s'))\\
&=1.
\end{align*}
However, for each $x<x_0$, $L(y,z)\leq x<x_0=L(r,s)$ implies that either $y<r$ or $z<s$. Thus, it holds that
$$\tau_{T}^{L}(\vep(r,1),\vep(s,1))(x)=\sup_{L(y,z)\leq x}T(\vep(r,1)(y),\vep(s,1)(z))=0.$$
Thus, $\tau_{T}^{L}(\vep(r,1),\vep(s,1))$ is not left continuous at the point $x_0\in]0,\infty[$, hence it is not a d.d.f., which is a contradiction. Therefore, the condition $(LCS)$ holds for the t-conorm $L$.

 (3)$\implies$(2): It is sufficient to show that $\tau_{T}^{L}(f,g)$ is left continuous on the open interval $]0,\infty[$. Given a point $x_0\in]0,\infty[$, for any $p<\tau_{T}^{L}(f,g)(x_0)=\sup_{L(r,s)\leq x_0}T(f(r),g(s))$, there are  real numbers $r_0$ and $s_0$ in $]0,\infty[$ such that $L(r_0,s_0)\leq x_0$ and $p<T(f(r_0),g(s_0))$. Since the function $T(f(-),g(-))$ is left continuous at the point $(r_0,s_0)$, there is some  $r_1$ in the open interval $]0,r_0[$ and some $s_1\in]0,s_0[$ such that   $p<T(f(r_1),g(s_1))$. Because   $L$ satisfies the condition $(LCS)$,   we have that 
 $$x=L(r_1,s_1)<L(r_0,s_0)\leq x_0.$$
 Thus, it holds that
 $$p<T(f(r_1),g(s_1))\leq\sup_{L(r,s)\leq x}T(f(r),g(s))=\tau_{T}^{L}(f,g)(x)\leq\tau_{T}^{L}(f,g)(x_0).$$
 Therefore, the map $\tau_{T}^{L}(f,g)$ is left continuous at any given point $x_0\in]0,\infty[$. 
\end{proof}

Based on Proposition \ref{tau_T^L vs tau_T,L} and Theorem \ref{tens vs tau}, we establish the following results, which refine and extend those presented in Proposition \ref{tauTL is tria func} and Proposition \ref{tauTL is tria func refo}.

\begin{thm}\label{tauTL} Let $T$ be a left continuous t-norm on $[0,1]$ and $L$ be a continuous t-conorm on $[0,\infty]$. Then the following statements are equivalent:
\begin{itemize}
\item[(1)] $\tau_{T,L}$ is a triangle function on $\Delp$;
\item[(2)] $\tau_{T,L}=L\otimes T$;
\item[(3)] $L$ satisfies the property (LCS).
\end{itemize}
\end{thm}

\section{Subsemigroups  of $(\Delp,L\otimes T)$}

 Recall from \cite{Schweizer1983} that a d.d.f. $f:[0,\infty]\lra[0,1]$ is \emph{non-defective} if $f(\infty^-)=1$, and it is \emph{strict} if it is continuous and strictly increasing on $[0,\infty]$.

In this section, we investigate several subsets of $\Delp$:
\begin{itemize}
\item $\CDp$, consisting of all non-defective d.d.f.'s.
\item $\CDp_0$, consisting of all d.d.f.'s in $\Delp$ which are continuous on $]0,\infty]$;
\item $\CDp_c$, consisting of all d.d.f.'s in $\CDp$ which are continuous on $[0,\infty]$. 
\item $\CDp_{sc}$, consisting of all strict d.d.f.'s in $\Delp$.
\end{itemize}
The hierarchical inclusion relations among the subsets are given by
$$\CDp_{sc}\subseteq \CDp_c\subseteq\CDp_0\subseteq\CDp\subseteq\Delp.$$
We demonstrate that the closedness  of these subsets under the operation $L\otimes T$ is fundamentally determined by the structural properties of the t-conorm $L$.

\begin{thm}\label{CDp clos} Let $T$ be a left continuous t-norm on $[0,1]$ and $L$ be a right continuous t-conorm on $[0,\infty]$. Then $L\otimes T(f,g)\in\CDp$ for all $f,g\in\CDp$ if and only if $L$ has no zero divisor. In this case, $(\CDp,L\otimes T)$ is a submonoid of $(\Delp,L\otimes T)$. 
\end{thm} 
\begin{proof} The ``only if'' part: If $L$ has a zero divisor $r_0\in]0,\infty[$, then there is some $s_0\in]0,\infty[$ such that $L(r_0,s_0)=\infty$. Let $x_0=\max\{r_0,s_0\}$, then $L(x_0,x_0)=\infty$. If $L(r,s)<x$ with $x<\infty$, then it holds that either $r\leq x_0$ or $s\leq x_0$. Take a function $f\in\CDp$ which is increasing strictly, then 
$$L\otimes T(f,f)(x)=\sup_{L(r,s)<x}T(f(r),f(s))\leq T(f(x_0),f(x_0))<1$$
for all $x<\infty$. Thus, $L\otimes T(f,f)(\infty^-)\leq T(f(x_0),f(x_0))<1$ and it is not in $\CDp$, a contradiction. Hence, $L$ has no zero divisor. 

The ``if'' part: We need check that $L\otimes T(f,g)(\infty^-)=1$ for all $f,g\in\CDp$. For any given $p<1$, there is some $q\in]p,1[$ such that $p<T(q,q)$ since $T$ is left continuous. Moreover, there is some $x_p\in]0,\infty[$ such that $q<f(x_p)$ and $q<g(x_p)$. Because $L$ has no zero divisor, we has that $L(x_p,x_p)<\infty$. Take a real number $y\in]L(x_p,x_p),\infty[$, then it holds that 
$$p<T(q,q)\leq T(f(x_p),g(x_p))\leq\sup_{L(r,s)<y}T(f(r),g(s))=L\otimes T(f,g)(y).$$
By the arbitrariness of $p$, we have that $L\otimes T(f,g)(\infty^-)=1$ as desired. 
\end{proof}

\begin{prop}\label{nece T} Let $T$ be a left continuous t-norm on $[0,1]$ and $L$ be a right continuous t-conorm on $[0,\infty]$.
 If $L\otimes T(f,g)\in\CDp_0$  for all  $f\in\CDp_0$ and all $g\in\CDp_c$, then   $T$ is continuous.
\end{prop}
\begin{proof}   For any give $p\in[0,1]$, let
$$f(x)=\begin{cases} 0& x=0,\\
                   p& 0<x\leq 1,\\
                   \min\{ p\cdot x, 1\} & 1<x\leq\infty, 
                   \end{cases}
$$
and let  
$$g(x)=\begin{cases}x& 0\leq x\leq 1,\\
                    1& 1<x\leq\infty.
       \end{cases}
$$
Clearly,   $f$ is in $\CDp_0$ and $g$ is in $\CDp_c$. 

Since $L$ is right continuous, it holds that $L(0^-,s)=L(0,s)=s$ for each $s\in [0,1]$. Thus, for any given $s$  and $x$ in $[0,1]$ with $s<x$,  there is some $r\in]0,x[$   such that $s\leq L(r,s)< x$. Hence, we can calculate that, for all $x\in]0,1]$,
$$L\otimes T(f,g)(x)=\sup_{L(r,s)<x}T(f(r),g(s))=\sup_{s<x}T(p,s)=T(p,x).$$

Moreover, $L\otimes{T}(f,g)(0^+)=T(0^+,p)=0$. Thus, if $L\otimes{T}(f,g)\in\CDp_0$, then   it is continuous on $[0,1]$. Therefore, $T(p,-)$ is   continuous   on $[0,1]$ for each $p\in[0,1]$. It follows that $T$ is a continuous t-norm.
\end{proof}

\begin{prop} \label{nece LS LCS}
 Let $T$ be a left continuous t-norm on $[0,1]$ and $L$ be a right continuous t-conorm on $[0,\infty]$.   
 If $L\otimes T(f,g)\in\CDp_0$ for all continuous d.d.f.'s $f$ and $g$, then $L$ satisfies the condition $(LS)$.
\end{prop}
\begin{proof}  
 If $L$  does not satisfy the condition $(LS)$, then   there are $r<r'$ and $s<s'$ in the  interval $]0,\infty[$ such that 
$$L(r,s)=L(r',s').$$ 
In this case, we show that there are  continuous d.d.f.'s $f$ and $g$ such that $h=L\otimes T(f,g)$ but it is not in $\CDp_0$. Our strategy is to show that  
$$h^\vee(y)=\inf_{T(p,q)>y}L(f^\vee(p),g^\vee(q)), 0\leq y\leq 1$$  
does not increase strictly on the interval $]h(0^+),1[$.  

Given any real number $y_0\in]0,1[$, define  continuous d.d.f.'s:  
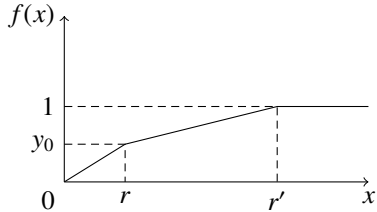
\begin{center}
\begin{tikzpicture}
        \draw [->](0,-1)--(4,-1) node[below]{$x$};
        \draw[->] (0,-1)--(0,1.2) node[left]{$f(x)$};
        \node[below left] at (0,-1){$0$};
        \node[below] at (0.8,-1){$r$};
        \node[below] at (2.8,-1){$r'$};
        \node[left] at (0,-0.5){$y_0$};
        \node[left] at (0,0){$1$};
        \draw (0,-1)--(0.8,-0.5)--(2.8,0)--(4,0);   
        \draw[densely dashed] (0,-0.5)--(0.8,-0.5)--(0.8,-1);
        \draw[densely dashed] (0,0)--(2.8,0)--(2.8,-1);
\node at(8,-.2)
{$f(x)=\begin{cases}\dfrac{y_0}{r}x & 0\leq x\leq r\\
                    y_0+\dfrac{1-y_0}{r'-r}(x-r) & r<x\leq r'\\
                    1 & r'\leq x\leq\infty
       \end{cases}
$};
\end{tikzpicture} 
\end{center}

and 
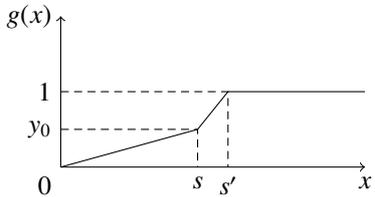
\begin{center}
 \begin{tikzpicture}
        \draw [->](0,0)--(4,0) node[below]{$x$};
        \draw[->] (0,0)--(0,2) node[left]{$g(x)$};
        \node[below left] at (0,0){$0$};
        \node[below] at (1.8,0){$s$};
        \node[below] at (2.2,0){$s'$};
        \node[left] at (0,0.5){$y_0$};
        \node[left] at (0,1){$1$};
        \draw (0,0)--(1.8,0.5)--(2.2,1)--(4,1);   
        \draw[densely dashed] (0,0.5)--(1.8,0.5)--(1.8,0);
        \draw[densely dashed] (0,1)--(2.2,1)--(2.2,0);
\node at(8,.8){$g(x)=\begin{cases}\dfrac{y_0}{s}x & 0\leq x\leq s\\
                    y_0+\dfrac{1-y_0}{s'-s}(x-s)& s<x\leq s'\\
                    1 & s'\leq x\leq\infty
       \end{cases}.$};
\end{tikzpicture}      
\end{center}

\begin{center}

\end{center}

Clearly, both $f$ and $g$ are continuous d.d.f.'s.  We can easily see the following facts:
\begin{enumerate}
\item[(a)]   $h(0^+)=0$, since  $h(x)\leq\min\{f(x),g(x)\}$  and $f(0^+)=\min\{g(0^+),f(0^+)\}=0$;  
\item[(b)] $f^\vee(y_0)=r$ and $f^\vee(1^-)=r'$;
\item[(c)] $g^\vee(y_0)=s$ and $g^\vee(1^-)=s'$.
\end{enumerate}

On one hand, if $T(p,q)>y_0$, then we have that both $p>y_0$ and $q>y_0$. Thus, it holds that 
\begin{align*}
h^\vee(y_0)
&=\inf_{T(p,q)>y_0}L(f^\vee(p),g^\vee(q))\\
&\geq L(f^\vee(y_0),g^\vee(y_0))\\
&=L(r,s).
\end{align*}
On the other hand, since $T$ is left continuous at the point $(1,1)$, for each $y\in]y_0,1[$,  there is some $p_y$ in $]y,1[$ such that 
$$y<T(p_y,p_y)<1.$$
 Thus, it follows that  
\begin{align*}
h^\vee(y)&=\inf_{T(p,q)>y}L(f^\vee(p),g^\vee(q))\\
&\leq L(f^\vee(p_y),g^\vee(p_y))\\
&\leq L(f^\vee(1^-),g^\vee(1^-))\\
&=L(r',s').
\end{align*}
So we have that 
$$L(r,s)\leq h^\vee(y_0)\leq h^\vee(y)\leq L(r',s'),$$ 
which means that $h^\vee(y_0)=h^\vee(y)$ because $L(r,s)=L(r',s')$. 
Therefore, the function $h^\vee$ is not strictly increasing on the interval $]0,1[$.  By Proposition \ref{cont vs stri incr}, the function $h=L\otimes T(f,g)$ is not continuous on the interval $(0,\infty]$, which is a contradiction.  


\end{proof}

\begin{lem}\label{lowe k} Let $T$ be a continuous t-norm on $[0,1]$ and $L$ be a right continuous t-conorm on $[0,\infty]$. Given functions $u$ and $v$ in $\Nabp$, define a function 
$$k(y)=\inf_{L(p,q)=y}L(u(p),v(q)), \quad y\in[0,1],$$  then for all $y'<y$ in $[0,1]$, there are real numbers $p'$ and $q'$ in $[0,1]$ such that \[y'<T(p',q')\quad\text{and}\quad
L(u(p'),v(q'))\leq k(y).\] 
\end{lem}
\begin{proof}
Since $k(y)$ is the infimum of the set $\{L(u(p),v(q))\}$, for each $n\in \mathbb{Z}^+$, there is some point $(p_n,q_n)$ in the set 
$$C_y=\{(p,q)\in[0,1]^2\mid T(p,q)=y\}$$ such that $$L(u(p_n),v(q_n))<k(y)+\frac{1}{n}.$$
Notice that the set $C_y$ is a closed subset of the square $[0,1]^2$ because $T$ is continuous, hence it is compact. So, the sequence $\{(p_n,q_n)\}_{n\in\mathbb{Z}^+}$ has a subsequence $\{(p_{n_k},q_{n_k})\}_{k\in\mathbb{Z}^+}$ converging  to a point  $(p_y,q_y)\in C_y$. Let $$p'_{n_k}=\inf\{p_{n_l}\mid l\geq k\},\quad q'_{n_k}=\inf\{q_{n_l}\mid l\geq k\}.$$  Clearly,  the sequence  $\{(p'_{n_k},q'_{n_k})\}_{k\in\mathbb{Z^+}}$ satisfies the properties in the below:
\begin{enumerate}
\item[(a)] both $\{p'_{n_k}\}_{k\in{\mathbb{Z}^+}}$ and $\{q'_{n_k}\}_{k\in\mathbb{Z^+}}$ are increasing;
\item[(b)] for all $k\in\mathbb{Z}^+$, $p'_{n_k}\leq p_{n_k}$ and $q'_{n_k}\leq q_{n_k}$;  
\item[(c)] the sequence $\{(p'_{n_k},q'_{n_k})\}_{k\in\mathbb{Z^+}}$
converges to $(p_y,q_y)$. 
\end{enumerate}
Thus, for all $k\in\mathbb{Z}^+$, it holds that $$p'_{n_k}\leq p_y, \quad q'_{n_k}\leq q_y$$ and $$\lim_{k\ra\infty}L(u(p'_{n_k}),v(q'_{n_k}))\leq\lim_{k\ra\infty}L(u(p_{n_k}),v(q_{n_k}))=k(y).$$ 
Since $y'<y$, by the continuity of $T$, there is some $(p'_{n_k},q'_{n_k})$ such that $$y'<T(p'_{n_k},q'_{n_k})\leq T(p_y,q_y)=y.$$ Moreover, we can easily see that   
$$L(u(p'_{n_k}),v(q'_{n_k}))\leq k(y).$$
\end{proof}

\begin{thm} \label{Delpc clos} Let $T$ be a left continuous t-norm on $[0,1]$ and $L$ be a right continuous t-conorm on $[0,\infty]$.  Then $L\otimes T(f,g)\in\CDp_0$  for all $f,g$ in $\CDp_0$ if and only if $T$ is continuous and $L$ satisfies the condition $(LS)$. In this case, $(\CDp_0,L\otimes T)$ is a submonoid of $(\Delp,L\otimes T)$. 
\end{thm}
\begin{proof} 
 The ``only if part'' holds by Proposition \ref{nece T} and Proposition \ref{nece LS LCS}. We need to check the ``if'' part.  Let $h=L\otimes T(f,g)$ for given $f,g\in\CDp_0$. To show that $h$ is also in $\CDp_0$, one need to show that $h^\vee$ is increasing strictly on the interval $]h(0^+),1[$.    
Define a function 
$$ k(y)=\inf_{T(p,q)=y}L(f^\vee(p),g^\vee(q)),$$
then $h^\vee=k^+$. Hence, $h^\vee$ is increasing strictly on $]h(0^+),1[$ if and only if $k$ also is. 
We claim that, for all $y'<y$ in $]h(0^+),1[$, it holds that $k(y')<k(y)$ indeed. We check it by 5 steps. 

\begin{enumerate}
\item[\emph{Step 1}:] By Lemma \ref{lowe k}, there is some point $(p',q')$ in the square $[0,1]^2$ such that $$y'<T(p',q')\leq y,\quad L(f^\vee(p'),g^\vee(q'))\leq k(y).$$ By the continuity of $T$, 
there is some $\xi\in]0,q'[$ such that $$T(p',\xi)=\dfrac{y'+T(p',q')}{2},$$ 
and   some $\eta\in]0,p'[$ such that 
$$T(\eta,\xi)=y'.$$

\item[\emph{Step 2}:] Notice that, on one hand, 
\begin{align*}
h(0^+)&=L\otimes T(f,g)(0^+)\\
&\geq L\otimes{T}(\vep(0,f(0^+)),\vep(0,g(0^+)))(0^+)\\
&=\vep(0,T(f(0^+),g(0^+)))(0^+)\\
&=T(f(0^+),g(0^+)),
\end{align*}
and on the other hand, 
\begin{align*}
h(0^+)&=\inf_{x>0}\sup_{L(r,s)<x}T(f(r),g(s))\\
&\leq\inf_{x>0} T(f(x),g(x))\\
&=T(f(0^+),g(0^+)).
\end{align*}
Thus, it holds that 
$$h(0^+)=T(f(0^+),g(0^+)).$$

\item[\emph{Step 3}:] Since $T(\eta,\xi)=y'>h(0^+)=T(f(0^+),g(0^+))$,  we have that  either $f(0^+)<\eta$ or $g(0^+)<\xi$. Without loss of generality, let $f(0^+)<\eta$, then $0<f^\vee(\eta)$. By $f^\vee$ being increasing strictly on $]f(0^+),1[$, it follows that
$$f^\vee(\eta)<f^\vee(p').$$ 

\item[\emph{Step 4}:] We check the  strict inequality $$L(f^\vee(\eta),g^\vee(\xi))<L(f^\vee(p'),g^\vee(q'))$$
in two cases.
\begin{enumerate}
\item[\emph{Case 1}:] If $g^\vee(\xi)<g^\vee(q')$, then by the condition $(LS)$, the strict inequality  holds as desired. 

\item[\emph{Case 2}:] If $g^\vee(\xi)= g^\vee(q')$, then by $g^\vee$ being   increasing strictly on $]g(0^+), \infty[$, it holds that $\xi<q' \leq g(0^+)$, hence $0=g^\vee(\xi)=g^\vee(q')$. In this case,   it follows that 
$$L(f^\vee(\eta),g^\vee(\xi))=f^\vee(\eta)<f^\vee(p')=L(f^\vee(p'),g^\vee(q')),$$
that is, the strict inequality also holds as desired.  
\end{enumerate}

\item[\emph{Step 5}:] We calculate that 
$$k(y')=\inf_{T(p,q)={y'}}L(f^\vee(p),g^\vee(q))
     \leq L(f^\vee(\eta),g^\vee(\xi))<L(f^\vee(p'),g^\vee(q'))\leq k(y).$$ 
\end{enumerate}
\end{proof}


\begin{thm} \label{cont ddf} Let $T$ be a  continuous t-norm on $[0,1]$ and $L$ be a right continuous t-conorm on $[0,\infty]$.  Then $L\otimes T(f,g)\in\CDp_c$   for all   $f,g\in\CDp_c$   if and only if  $L$ satisfies the condition $(LS)$. In this case, $(\CDp_c,L\otimes T)$ is a subsemigroup, but not a submonoid of $(\Delp,L\otimes T)$ since the identity $\vep(0,1)$ is not in $\CDp_c$.
\end{thm}
\begin{proof}
The ``only if'' part is from Proposition \ref{nece LS LCS}. To show the ``if'' part, since the condition $(LS)$ is satisfied, by Theorem \ref{Delpc clos}, $L\otimes T(f,g)$ is already in $\CDp_0$. Additionally, it always holds that 
$L\otimes T(f,g)\leq f$. Hence, $L\otimes T(f,g)(0^+)=f(0^+)=0$, which means that $L\otimes T(f,g)$ is continuous at $0$. Therefore, $L\otimes T(f,g)$ is also continuous if both $f$ and $g$ are.
\end{proof}

\begin{rem}In the above theorem, we assume that the t-norm $T$ is continuous, but whether it is necessary remains open. 
\end{rem}

\begin{cor}\label{continc ddf}
    Let $T$ be a continuous t-norm on $[0,1]$ and $L$ be a continuous t-conorm on $[0,\infty]$. Then  $L\otimes T(f,g)\in\CDp_{sc}$   for all   $f,g\in\CDp_{sc}$  if and only if   $L$ satisfies the condition $(LS)$ and $T$ satisfies the condition  
    \[x<y,x'<y'\Longrightarrow T(x,y)<T(x',y').\leqno{(TS)}\] In this case, $(\CDp_{sc},L\otimes T)$ is a subsemigroup of $(\Delp,L\otimes T)$.
\end{cor}
\begin{proof}
    Define a binary operation $T\otimes{L}:\Nabp\times{\Nabp}\lra{\Nabp}$ by $$T\otimes{L}(f^\vee,g^\vee)(y)=(L\otimes{T}(f,g))^\vee(y)=\inf_{T(p,q)>y}L(f^\vee(p),g^\vee(q))$$ for all $f,g\in\Delp$ and $y\in[0,1].$ By the dual of Theorem \ref{cont ddf}, when $L$ is continuous, $T\otimes{L}(f^\vee,g^\vee)$ is continuous for all continuous $f^\vee,g^\vee$ in $\Nabp$ if and only if $T$ satisfies the condition $(TS)$, which is exactly equivalent to that $L\otimes{T}(f,g)$ is strictly increasing for all strictly increasing d.d.f.'s $f,g\in\Delp$.
\end{proof}

\begin{thm}\label{CDpc clos} Let $T$ be a left continuous t-norm on $[0,1]$ and $L$ be a right continuous t-conorm on $[0,\infty]$.  Then $L\otimes T(f,g)\in\CDp_c$  for all   $f\in\CDp$ and all $g\in\CDp_c$ if and only if $T$ is continuous and $L$ satisfies the  cancellation law. In this case, $(\CDp_c,L\otimes T)$ forms an ideal of $(\CDp,L\otimes T)$. 
\end{thm}
\begin{proof} 
The  ``only if'' part:  Firstly, by Proposition \ref{nece T}, the t-norm $T$ is continuous. Secondly, we show that $L$ satisfies the cancellation law. Otherwise,  there is some $r\in]0,\infty[$ and $s<s'$ in the  interval $]0,\infty[$ such that 
$$L(r,s)=L(r,s').$$ 
In this case, we show that there is a  function $f$ in $\CDp$ and a  function $g$ in $\CDp_c$  such that the function $h=L\otimes T(f,g)$ is not continuous, 
that is, the function 
$$h^\vee(y)=\inf_{T(p,q)>y}L(f^\vee(p),g^\vee(q)), 0\leq y\leq 1$$   does not increase strictly on the interval $[0, 1]$.  

Fix a real number  $y_0$ in $]0,1[$. 
Choose a  continuous  functions $g$ in $\Delp$ such that $g$ is increasing strictly on $[0,s']$ and $g(s)=y_0$, $g(s')=1$, and choose   another function $f$ in $\Delp$ such that $f(r)=y_0$ and $f(r^+)=1$. 
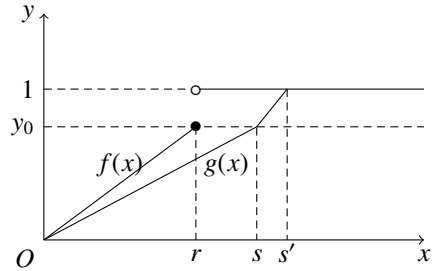
\begin{center}
\begin{tikzpicture}
\draw [->](-1,0)--(4,0) node[below]{$x$};
\draw[->] (-1,0)--(-1,3) node[left]{$y$};
\node[below left] at (-1,0){$O$};
\node[below] at (1.8,0){$s$};
\node[below] at (2.2,0.08){$s'$};
\node at(1.4,1) {$g(x)$};
\node[left] at (-1,1.5){$y_0$};
\node[left] at (-1,2){$1$};
\draw (-1,0)--(1.8,1.5)--(2.2,2)--(4,2);   
\draw[densely dashed] (-1,1.5)--(1.8,1.5)--(1.8,0);
\draw[densely dashed] (2.2,2)--(2.2,0);
\draw[densely dashed] (1.8,1.5)--(4,1.5);
\draw(-1,0)--(1,1.5) node{$\bullet$};
\node at(1,1.98){$\circ$};
\draw[densely dashed] (1,1.5)--(1,0) node[below]{$r$};
\draw[densely dashed] (-1,2)--(.9,2);
\draw(1.05,2)--(2.2,2);
\node at(0,1){$f(x)$};
\end{tikzpicture}
\end{center}
Then we see that
$$f^\vee(y_0)=f^\vee(1^-)=r,\quad  g^\vee(y_0)=s,\quad g(1^-)=s'.$$ 
We consider the values $h^\vee(y_0)$ and $h^\vee(y)$ for all $y\in]y_0,1[$. Notice that for all $p,q\in[0,1]$, 
if $t<T(p,q)<1$, then
$$t<T(p,q)\leq\min\{p,q\}<1.$$
Thus, it follows  that 
$$L(r,s)=L(f^\vee(y_0),g^\vee(y_0))
\leq\inf_{T(p,q)>y_0}L(f^\vee(p),g^\vee(q))=h(y_0),$$
and  
$$h(y_0)\leq h(y)=\inf_{T(p,q)>y}L(f^\vee(p),g^\vee(q))\leq L(f^\vee(1^-),g^\vee(1^-))=L(r,s').$$
So, we have that 
$$L(r,s)=h^\vee(y_0)=h^\vee(y)=L(r,s')$$ for all $y$ in $]y_0,1[$. 
Therefore,  the function $h^\vee$ is not strictly increasing. That means, $h=L\otimes T(f,g)$ is not continuous, which is a contradiction.  

The ``if'' part: We need to show that $(L\otimes T(f,g))^\vee$ is increasing strictly on $[0,1]$.  Firstly, we verify that the function 
$$h(y)=\inf_{T(p,q)=y}L(f^\vee(p),g^\vee(q))$$
is increasing strictly on $[0,1]$. Let $y'<y$ in $[0,1]$, by Lemma \ref{lowe k}, there are real numbers $p'$ and $q'$ in $[0,1]$ such that 
$$y'<T(p',q')\leq y,\quad L(f^\vee(p'),g^\vee(q'))\leq h(y).$$
By the continuity of $T$, 
there is some $\xi\in]0,q'[$ such that $$T(p',\xi)=\dfrac{y'+T(p',q')}{2},$$ 
and   some $\eta\in]0,p'[$ such that 
$$T(\eta,\xi)=y'.$$
 
Since $f(\infty^-)=1$, one has that $f^\vee(\eta)<\infty$. When $L$ satisfies the cancellation law, one has that 
$$L(f^\vee(\eta),g^\vee(\xi))<L(f^\vee(\eta),g^\vee(q'))\leq L(f^\vee(p'),g^\vee(q'))$$ 
since the function $g^\vee$ is strictly increasing. Thus, we obtain  that
\begin{align*} h(y')&=\inf_{T(p,q)={y'}}L(f^\vee(p),g^\vee(q))\\
     &\leq L(f^\vee(\eta),g^\vee(\xi))\\
     &< L(f^\vee(p'),g^\vee(q'))\\
      &\leq h(y).
 \end{align*}
 That is, the function $h$ is increasing strictly on $[0,1]$ as desired. Secondly, since  $(L\otimes T(f,g))^\vee=h^+$, the function $(L\otimes T(f,g))^\vee$ is also increasing strictly on $[0,1]$, which implies that $L\otimes T(f,g)$ is continuous.  
\end{proof}

In the above theorem, the condition $f\in\CDp$ is necessary. Otherwise, let $f=\vep(0,p)$ with $p\in[0,1[$, then for all $f\in\Delp$, 
$$L\otimes T(f,g)(x)=\sup_{L(r,s)<x}T(f(r),g(s))\leq p$$ for all $x<\infty$. Thus, $L\otimes T(f,g)(\infty^-)\leq p<1$ and it is not continuous at $\infty$.

\section{Conclusion}
In this paper, we investigate triangle functions of the form $L\otimes T$,  which is the tensor product of a left continuous t-norm $T$ on $[0,1]$ and a right continuous t-conorm $L$ on $[0,\infty]$. The main results include that
\begin{enumerate}
\item[(1)] If $L$ is a continuous t-conorm on $[0,\infty]$, $L\otimes T=\tau_{T,L}$ $\iff$ $\tau_{T,L}$ is a triangle function $\iff$  $L$ satisfies the property $(LCS)$ (see Theorem \ref{tauTL}); 
\item[(2)]$(\CDp,L\otimes T)$ is a submonoid of $(\Delp,L\otimes T)$ $\iff$ $L$ has no zero divisor (see Theorem \ref{CDp clos}), which improves the result of Theorem 12.1 (b) in \cite{Saminger2008};
\item[(3)] $(\CDp_0,L\otimes T)$ is a submonoid of $(\Delp,L\otimes T)$ $\iff$ $T$ is continuous and $L$ satisfies the condition $(LS)$ (see Theorem \ref{Delpc clos});
\item[(4)] if $T$ is continuous, then 
 $(\CDp_c,L\otimes T)$ is a subsemigroup of $(\Delp,L\otimes T)$ $\iff$ $L$ satisfies the condition $(LS)$ (see Theorem \ref{cont ddf}); 
 \item[(5)] if $T$ and $L$ both are continuous, then 
 $(\CDp_{sc},L\otimes T)$ is a subsemigroup of $(\Delp,L\otimes T)$ $\iff$ $T$ satisfies the condition $(TS)$ and $L$ satisfies the condition $(LS)$ (see Corollary \ref{continc ddf});
\item[(6)]  $(\CDp_c,L\otimes T)$ is an ideal of $(\CDp,L\otimes T)$ $\iff$ $T$ is continuous and $L$ satisfies the cancellation law (see Theorem \ref{CDpc clos}). 
\end{enumerate} 

It is worth noting that Theorem \ref{CDpc clos} 
reveals a flaw in Corollary 7.5 of \cite{Saminger2008}. Specifically, taking $T=M$ and $L=M^*$, where $M^*$ satisfies the conditions stated in the corollary, we find that
$$L\otimes T(f,g)(x)=\min\{f(x),g(x)\},$$
which implies that $L\otimes T(f,g)=f$ is discontinuous when $f$ is a discontinuous distance distribution function lying below a continuous distance distribution function $g$.

\section*{Declaration of competing interest}
The authors declare that they have no known competing financial interests or personal relationships that could have appeared to influence the work reported in this paper.
\section*{Acknowledgment}
The authors acknowledge the support of National Natural Science Foundation of China (12171342). 

The authors are grateful to the anonymous referees for their helpful remarks and suggestions.

\end{document}